\documentclass[11pt]{amsart}
\usepackage{amscd,amssymb}
\theoremstyle{plain}

\newtheorem{Thm}[equation]{Theorem}

\newtheorem{Prop}[equation]{Proposition}

\numberwithin{equation}{section}

\newcommand{\mc}[1]{{}}

\newcommand{\q}{\mathbb{Q}}

\renewcommand{\q}{\mathbb{Q}}
\newcommand{\n}{\mathbb{N}}

\newcommand{\br}{\mathbb{R}}
\newcommand{\qp}{\mathbb{Q}_p}

\newcommand{\ba}{\backslash}

\newcommand{\G}{\Gamma}

\newcommand{\CG}{\operatorname{Comm}(\Gamma )}
\newcommand{\Cal}{\mathcal}
\renewcommand{\P}{\mathcal P}

\renewcommand{\gg}{\Gamma\ba G}

\newcommand{\supp}{\operatorname{supp}}

\newcommand{\Comm}{\operatorname{Comm}}

\renewcommand{\deg}{\text{deg}}

\begin{document}

\title[Equidistribution of Hecke points]{Ergodic theoretic
proof of equidistribution of Hecke points}
\author{Alex Eskin and Hee Oh}
\address{Mathematics Department, University of Chicago, Chicago,
IL 60637}

\thanks{The first author is partially supported by Packard foundation.}
\thanks{The second author partially supported by NSF grants DMS 0070544
and DMS 0333397.}

\email{eskin@math.uchicago.edu}
\address{Mathematics Department, Princeton University, Princeton,
NJ 08544, Current address: Math 253-37, Caltech, Pasadena, CA 91125}
\email{heeoh@its.caltech.edu}
\maketitle

\section{Introduction}
Let $G$ be a connected non-compact $\q$-simple real algebraic group
defined over $\q$, that is, the identity component of the group of
the real points of a connected $\q$-simple algebraic group which is
$\br$-isotropic.
Let $\G\subset G(\q) $ be an arithmetic subgroup of $G$.
As is well known, $\G$ has finite co-volume in $G$ [BH].
 Denote by $\mu_G$ the $G$-invariant Borel probability
measure on $\gg$.
Two subgroups $\G_1$ and $\G_2$ of $G$ are said to be
{\it commensurable} with each other if $\G_1\cap \G_2$ has a finite
index both in $\G_1$ and $\G_2$.
 The commensurator group $\CG$ of $\G$ is defined as follows:
$$\CG=\{ g\in G : \G \text{ and } g\G g^{-1} \text{ are commensurable
with each other} \} .$$
Since $\G$ is an arithmetic subgroup, $\CG$ contains $G(\q)$
and in particular,
$\G$ has an infinite index in $\CG$.

For an element $a\in \text{Comm}(\G)$,
the $\G$-orbit $\G\ba \G a\G$ in $\G\ba G$ has finitely many
points called {\it Hecke points} associated with $a$.
We set
 $$\deg(a)=\# \G\ba \G a\G .$$
 It is easy to see that $\deg(a)=[\G:\G\cap a^{-1}\G a]$.

In this paper, we are interested in the
equidistribution problem of the Hecke points $\G\ba \G a\G$
as $\deg(a) \to \infty$. Namely,
for any continuous function $f$ in $\G\ba G$ with compact support,
 any $x\in \G\ba G$, and for any sequence $a_i\in \Comm (\G)$
with $\deg(a_i)\to \infty$,

\begin{equation}\label{h}\text{Does }\operatorname{T}_{a_i}(f)(x):=
\frac{1}{\deg(a_i)}\sum_{\gamma \in \G\ba \G a_i\G}f(\gamma x)
\text{ converge to }
\int_{\G\ba G} f \, d\mu_G ?\end{equation}

We remark that $\operatorname{T}_{a_i}(f)$ is well defined as a function on $\G\ba G$.
When $G$ is simple, $a_i\in G(\q)$ and $\G$ is a {\it congruence} subgroup,
this was answered in the affirmative way in [COU],
based on the adelic interpretation of $L^2(\gg)$ and an information on
the local harmonic analysis of the $p$-adic groups $G(\qp)$.
For smooth functions, the methods in [COU] also give rate of the convergence.
Several partial results in this direction
 were known (see [C], [CU], [Sa], [GM], etc.)

It was pointed out by Burger and Sarnak [BS, "Theorem 5.2"] in 1991
that the equidistribution of Hecke points follows from Ratner's
measure classification theorem [Ra] provided $a_i$'s
 converge to an element not belonging to
$\CG$.  However Burger and Sarnak did not give a detailed proof of
this claim: this was done by Dani and Margulis [DM2, Corollary 6.2]
in 1993, where they deduced the assertion from their ergodic results
built up on Ratner's fore-mentioned theorem.

Unlike the method used in [COU], this ergodic theoretic method
does not provide a rate information of the equidistribution.
 However it works
 for much more general cases, for instance, $\G$ can be a
 non-congruence subgroup and $a_i$'s need not be necessarily
 elements of $G(\q)$.

Our main purpose of this paper is to present an ergodic theoretic
proof of the following result:
\begin{Thm}\label{T} Let $G$ be a connected non-compact
 $\q$-simple real algebraic group
defined over $\q$ and $\G \subset G(\q)$ an
arithmetic subgroup of $G$.
 Let $\{a_i \in \CG\}$ be a sequence such that
$\lim_{i \to \infty} \operatorname{deg}(a_i)=\infty $. Then for
any bounded continuous function $f$ on $\gg$ and for any $x\in \gg$,
$$\lim_{i \to \infty}
\operatorname{T}_{a_i}(f)(x) = \int_{\gg} f (g) \,d\mu_{G} (g). $$
\end{Thm}

\section{Limits of $H$-invariant measures}
 Let $G$ be a connected real semisimple
 algebraic group defined over $\q$ and
$\G\subset G(\q)$ an arithmetic subgroup of $G$. Let $H$ be a
connected real non-compact semisimple $\q$-simple
subgroup of $G$. Then $\G\cap H$ is an irreducible ($H$ being
$\q$-simple) Zariski dense lattice in $H$.

 Let $\{g_m \in G\}$ be a sequence such that
$g_m^{-1} \G g_m \cap  H$ is commensurable with $\G\cap H$.
It follows that
each $\G\ba \G g_mH$ is closed [Rag] and there exists
  the unique
 $H$-invariant probability
 measure, say $\nu_m$, in $\G\ba G$ supported on $\G\ba \G g_m H$.

Let $Y$ denote $\G\ba G$ if $\G\ba G$ is compact; and otherwise
the one point compactification $\G\ba G \cup \{\infty\}$. The space $ \P
(Y)$ of the probability measures on
$Y$ equipped with the weak$^*$-topology is weak$^*$
compact.

Our basic tools for the proof of Theorem \ref{T}
 are the following two propositions:
Denote by $H_N$ the unique maximal connected normal subgroup of
$H$ with no compact factors.
\begin{Prop}\label{Non-div}
Suppose that
$g_mH_Ng_m^{-1}$ is not contained in
 any proper parabolic $\q$-subgroup
 of $G$ for each $m$.
 Then every weak
limit of $\{\nu_{m} :m\in \n\}$ in $\P (Y)$ is supported in $\G\ba G$.
\end{Prop}

This proposition is
 shown in [EO, Proposition 3.4] based on theorems of Dani and
Margulis ([DM1, Theorem 2] and [DM2, Theorem 6.1]).

\begin {Prop}\label{limit} Suppose that $\nu _m $ weakly
converges to a measure $\nu$ in $\P(\G\ba G)$ as $m \to \infty$. Then
there exists a closed connected subgroup $L$ of $G$ containing $H$
such that
\begin{itemize}
\item[(1)] $\nu$ is an $L$-invariant measure supported on $\G\ba \G c_0 L$ for
some $c_0\in G$;
\item[(2)] $\G \cap c_0Lc_0^{-1}$ is a Zariski dense lattice in $c_0Lc_0^{-1}$
and hence in particular $c_0Lc_0^{-1}$ is defined over $\q$;
\item[(3)] there exist $m_0\in \n$ and a sequence
$\{ x_m \in\Gamma g_m H\}$ converging to $c_0$ as $m \to \infty$
such that  $c_0Lc_0^{-1}$ contains the subgroup generated by $\{
x_m H x_m^{-1} : m \ge m_0\}$.
\end{itemize}\end{Prop}
This proposition is deduced from
the following theorem of Mozes and Shah:
\begin{Thm} [MS, Theorem 1.1]\label{MS}
Let $\{ u_i(t)\}_{t\in \br}$, $ i\in \n$
  be a sequence of unipotent one-parameter
subgroups of $G$ and let $\{ \nu_m: m\in \n\}$ be a sequence in
$\P (\G\ba G)$ such that each $\nu_i$ is an ergodic
$\{u_i(t)\}$-invariant measure. Suppose that $\nu_m \to \nu$ in
$\P (\G\ba G)$ and let $x\in \supp(\nu)$. Then the following holds:
\begin{itemize}
\item[(1)] $\supp(\nu)=x\Lambda (\nu) $ where
$\Lambda (\nu)=\{g\in G : \nu g=\nu\}$.

\item[(2)] Let $g_i' \to e$ be a sequence in $G$ such that
for every $i\in \n$, $xg_i'\in \supp(\nu_i)$ and the trajectory
$\{xg_i'u_i(t)\}$ is uniformly distributed with respect to
 $\nu_i$. Then there exists an $i_0\in \n$ such that for all $i \ge i_0$,
$$\supp(\nu_i)\subset\supp(\nu)g_i'.$$
\item[(3)] $\nu$ is invariant and ergodic for the action of
the subgroup generated by the set $\{g_i'u_i(t){g_i'}^{-1}:i\ge
i_0\}$.
\end{itemize}
\end{Thm}

\noindent{\bf Proof of Proposition \ref{limit}}
Since $g_m^{-1}\G g_m\cap H$ is an irreducible
lattice in $H$, every non-compact simple factor
of $H$ acts ergodically on each $\G\ba \G g_m H$ with
respect to $\nu_m$.
There exists a unipotent one-parameter subgroup $U:=\{u(t)\}$ in
$H_N$ not contained in any proper closed normal subgroup of
 $H_N$ (cf. Lemma 2.3 [MS]).
 Then by Moore's ergodicity theorem (cf. Theorem 2.1 in [BM]),
$U$ acts ergodically with respect to each $\nu_m$. Moreover by the
Birkhoff ergodic theorem, the following subset $R$ has the zero
co-measure in $H$:$$\{ h\in H: \G\ba \Gamma g_m hu(t) \text{ is
uniformly distributed in $\G\ba \Gamma g_mH$ w. r. t. $\nu_m$ for
each $m\in \n$}\}.$$ Hence for any $h\in R$ and for any continuous
bounded function $f$ on $X$ with compact support, we have
$$\lim_{T\to\infty}\frac{1}{T}\int_{0}^Tf(g_mhu(t))\, dt=\int_Xf\, d\nu_m .$$
If we set
 $L=\Lambda (\nu)$,
we have that $\nu$ is supported on $\G\ba \G c_0 L$ for some $c_0\in
\supp(\nu)$ by Theorem \ref{MS}(1).
 There exist $\gamma_m \in \G$ and
$h_m\in R$ such that $\gamma_m g_m h_m \to c_0$ as $ m \to
\infty$. If we set
$$x_m:=\gamma_m g_m h_m\quad\text{
 and }\quad g_m':=c_0^{-1}x_m ,$$
then $$\lim_{m\to \infty}g_m' = e$$ and
$\G\ba \Gamma c_0g_m'u(t)=\G\ba \Gamma
g_mh_m u(t)$ is uniformly distributed with respect to $\nu_m$. By
Theorem \ref{MS}(2), there exists an $m_0\in \n$ such that for all
$m\ge m_0$
 $$\supp (\nu_m)\subset \supp(\nu)g_m',\quad
\text{or equivalently }\quad \G\ba \G g_m H\subset \G\ba \G c_0Lg_m' .$$ Hence
$\Gamma x_mHx_m^{-1}\subset \Gamma (c_0Lc_0^{-1})$. By the
connectedness of $H$,
 we may assume that $c_0Lc_0^{-1}$ is connected and
$$  \{x_mHx_m^{-1}: m\ge m_0\} \subset c_0Lc_0^{-1}.$$
By Theorem \ref{MS}(3), the subgroup generated by the set $\{
x_mHx_m^{-1}: m\ge m_0\}$ acts ergodically on $\G c_0Lc_0^{-1}$.
Hence $c_0Lc_0^{-1}$ is the smallest closed subgroup containing
the subgroup generated by the set $\{ x_mHx_m^{-1}: m\ge m_0\}$
such that the orbit $\G c_0Lc_0^{-1}$ is closed. This proves (3).
The second claim (2) follows
 from [MS, Proposition 2.1].

\section{Proof of Theorem \ref{T}}
For a $G$-space $X$ and a subgroup $M$ of $G$, $\Cal P (X) ^M$
denotes the space of $M$-invariant Borel probability measures on
$X$. We recall the ergodic theoretic approach suggested in [BS].
Let $\Delta (G)$ be the diagonal embedding of $G$ into $G \times G$, that
is, $\Delta (G)=\{(g, g): g\in G\}$. For each $\nu\in \Cal P(\gg)^{\G}$,
the measure $\tilde \nu$ defined by
$$\tilde \nu(f):=\int _{\gg}\int_{\gg}
f(g, hg)\, d\mu _G (g) d\nu (h)\quad $$
 for any bounded continuous function $f$ on $\gg \times \gg$
is a $\Delta (G)$-invariant probability measure on $\gg \times \gg  $.
Moreover the map $\nu \mapsto \tilde\nu$ is a homeomorphism from
$\Cal P (\gg) ^{\G} $ to $ \Cal P (\gg \times \gg) ^{\Delta (G)} $. We make
the following simple observation:

\begin{itemize}
\item
For $a\in \CG$, if we set $\nu_{a}= \frac{1}{\deg(a)}\sum_{y\in
\G\ba \G a\G}\delta_{y}$ where $\delta_{y}$ denotes the standard
delta measure with support $y$, then $\tilde \nu_{a}$ is the
(unique) $\Delta (G)$-invariant probability measure supported on
$[(e,a)]\Delta (G)\subset (\G\times\G)\ba (G\times G)$.
\item The element $a\in G$ is contained in $\CG$ if and only if
the orbit $[(e, a)]\Delta (G)$ is closed and supports a finite
 $\Delta (G)$-invariant measure.
\end{itemize}

Set $X=(\G \times \G)\ba (G\times G)$ and consider its one point
compactification $X\cup \{\infty\}$. By the fact that $\P (X\cup
\{\infty\})$ is compact with respect to weak $^*$-topology and the
above observation,
 it suffices to show that, assuming
  the sequence $\{\tilde \nu_{a_i}\}$ weakly converging to $\tilde \nu$ in
$\P (X\cup\{\infty\})$, the limit $\tilde \nu$ is $G\times
G$-invariant and supported on $X$. Note that for each $i$, $(e,
a_{i}) \Delta (G) (e, a_i^{-1})\cap( \G \times \G)$ is commensurable with
$\Delta (G)\cap (\G \times \G)$. It follows that $(e,
a_{i}) \Delta (G) (e, a_i^{-1})$ contains a Zariski dense subset
contained in $\Delta (G)(\q)$. This implies that $(e, a_{i}) \Delta (G)(e, a_i^{-1})$ is
a $\q$-subgroup of $G\times G$ (cf. [Zi, Prop. 3.18]) and $\q$-simple
as well.
Moreover,
it is easy to see that the unique maximal
connected normal subgroup of $(e, a_{i}) \Delta (G) (e, a_i^{-1})$
with no compact factors cannot be contained in any proper parabolic
$\q$-subgroup of $G\times G$ for each $i$. Therefore by
Proposition \ref{Non-div}, $\tilde \nu$
 is supported on $X$.
Also by applying Proposition \ref{limit}, $\tilde \nu$ is either
$G\times G$-invariant or $\Delta (G)$-invariant
 supported on $x\Delta (G)$ for some $x\in X$.
Suppose that the latter case happens. Then Proposition \ref{limit}
 also says that
there exist $i$ and $y_j\in \G\times \G$ such that
$y_j(e, a_{j}) \Delta (G) (e, a_j^{-1})y_j^{-1}= (e, a_{i})
\Delta (G) (e, a_{i}^{-1})$ for all $j\ge i$.
That is, $(e, a_i^{-1})y_j (e, a_j) $ belongs to the normalizer
of $\Delta (G)$ in $G\times G$ for all $j\ge i$.
 Since $\Delta (G)$ has a finite index in
its normalizer in $G\times G$, we have
$(e, a_i^{-1})y_j (e, a_j) \in \Delta (G)$ and hence
 $[(e,a_{i})]\Delta (G)=[(e,a_{j})]\Delta (G)$ for infinitely many $j$.
This implies that $\G a_i =\G a_j$ and hence
$\deg(a_j)$ is constant for infinitely many $j$.
This is a
contradiction, since $\text{deg}(a_{i})$ tends to $\infty$ as $i
\to \infty$. Therefore $\tilde \nu$ is
 the $G\times G$-invariant
probability measure supported on $X$, as desired.

\enddocument